\documentclass[final]{siamltex}

\usepackage{graphics}
\usepackage{epsfig}
\usepackage[ps2pdf,hyperindex]{hyperref}
\usepackage{bbm}
\usepackage[latin1]{inputenc}
\usepackage{amsmath,amssymb}

\title{Alternatives for optimization in systems and control: convex and non-convex approaches
\thanks{This research is supported by the Belgian Network DYSCO
(Dynamical Systems, Control, and Optimization) funded by
the Interuniversity Attraction Poles Programme of the Belgian
State, Science Policy Office.}}

\author{Emile Simon\thanks{Applied Mathematics Department, 
        Université Catholique de Louvain, 4, avenue Georges Lemaître,
        1348, Louvain-la-Neuve, Belgium ({\tt Emile.Simon@uclouvain.be}).}}

\begin{document}

\maketitle



\pagestyle{myheadings}
\thispagestyle{plain}

In this presentation, we will develop a short overview of main trends of optimization in systems and control, and from there outline some new perspectives emerging today. More specifically, we will focus on the current situation, where it is clear that convex and Linear Matrix Inequality (LMI) methods have become the most common option. However, because of its vast success, the convex approach is often the only direction considered, despite the underlying problem is non-convex and that other optimization methods specifically equipped to handle such problems should have been used instead. We will present key points on this topic, and as a side result we will propose a method to produce a virtually infinite number of papers.

\vspace{0.2cm}
\hrule
\vspace{0.2cm}

After the origins of optimal control in the late fifties, championed by the works of Pontryagin and Bellman, the general intent was to try to find (almost) closed-form expressions for the considered problems. This perspective was pursued within the frameworks of Ricatti/Lyapunov equations approaches, with for instance the famous solution of Kalman for the LQ(R) problem. This period can be referred to as that of state-space methods and of so-called `modern control'. A summary about this trend can be found in the introduction of \cite{ZDG96} (see also ref. therein), and main results of that research direction are summarized for instance in \cite{AM89,B91}.

Somewhere in the eighties, the interest shifted to the more general framework of problems that can be formulated under the standard form of (nonlinear continuous) optimization problems:
\begin{equation}
\begin{split}
 & \min_x f(x)  \\
s.t.\, & x \in \Omega \subseteq \mathbb{R}^n
\end{split}
\tag{1}
\end{equation}
Such problems do not admit in general an analytical solution and must be solved through iterative methods. Problems of the type of (1) can be split between convex problems, where $f(x)$ and $\Omega$ are convex, and non-convex problems, where $f(x)$ and/or $\Omega$ are non-convex. The field of nonlinear programming, concerned with solving (1), is made of both these subsets. The classical reference book of that field is \cite{B99}, an alternative being \cite{LY08}.
\subsection*{The success of convexity}
We recommend the first chapter of \cite{SW05} for a clear and summarized but nevertheless quite comprehensive presentation of the main concepts and results of convex optimization. More detailed and formal presentations are developed in the following three fundamental books \cite{NN94,BN01,BV04}, that are the main references on convex programming (see also \cite{B99,LY08}, including convex programming).

\newpage

The great strength and main interest of convex optimization is that:
\vspace{0.1cm}

\textit{All locally optimal solutions of a convex problem are also globally optimal}.
\vspace{0.1cm}
 
This can be easily proven by observing that all locally optimal solutions must belong to the same convex set with the same objective value, which are then also globally optimal (see \cite[Subsec. 4.2.3]{BV04}). Moreover, most convex problems can be solved in polynomial time (cfr. e.g. \cite{NN94}).

The main class of convex problems is those that are classified as Semi-Definite Programs (SDPs). These problems can be defined by a linear objective function and by constraints admitting a Linear Matrix Inequality (LMI) representation:
\begin{equation}
\begin{split}
 & \min_x c^T x \\
s.t.\, & F(x) \preceq 0,
\end{split}
\tag{2}
\end{equation}

where $c \in \mathbb{R}^n$ and the LMI is defined as $F(x)=F_0+ \sum^n_{i=1} F_i x_i$, where all matrices $F$ belong to $\mathbb{S}^k$, the set of symmetric $k\times k$ matrices ($F(x) \preceq 0$ means that the maximum real part of the eigenvalues of $F(x)$ must be negative or zero).

This class of problems is now the most often considered when trying to solve a problem of optimization in systems and control.

There exist some historical introductions presenting the apparition and development of LMIs in systems and control. The most comprehensive is that given in the three pages of \cite[Sec. 1.2.]{BGFB94} (see also \cite[Subsec. 1.4.6]{SW05}). We will not copy this here and recommend these pages to the reader. The last crucial result of this development are the interior point methods extended by Nesterov and Nemirovsky to SDPs, in a most effective fashion (see \cite{NN94} and \cite{BN01}, the initial result was presented in 1988 in \cite{NN88}). Since then, many efficient solvers are available to solve SDPs (see Sedumi, LMILab, SDPT3 and other algorithms listed for instance in \cite{APH02}).

The availability of such solvers had the consequence in control of, as stated in \cite{SW05}: ``\textit{a significant paradigm shift:  Instead of arriving at an analytical solution of an optimal control problem and implementing such a solution in software so as to synthesize optimal controllers, the intention is to reformulate a given control problem to verifying whether a specific linear matrix inequality is solvable or, alternatively, to optimizing functionals over linear matrix inequality constraints.}''. The article \cite{VB97} is also recommended for more observations on this line.

The trend of the development of LMIs in control is still booming today. This explosion can be quantified for instance with the number of articles with the keyword LMI in the IEEEXplore database\footnote{See \url{http://www.eeci-institute.eu/pdf/M016/lmiI2.pdf}, p.9 for a graph.}: This explosion began around 1991-1994 (between the books \cite{BB91,BGFB94}), the 1st January 2009 there were 3133 cumulative hits and around 5300 hits at the beginning of April 2012.

Most of the main classical results of LMI formulations of systems and control optimization problems are presented in \cite{BB91,BGFB94,SGC97,EN00,SW05}.
\subsection*{Non-convex problems and methods}
An undesirable drawback of the success and large domination of convex optimization in systems and control is that it regularly happens that this direction is the only option considered, even though the underlying problem is non-convex (as such, or after de-approximating it).

A classical situation is that of trying to solve non-convex problems with iterative LMI algorithms (ILMIs), which could also be labelled as Sequential SDPs (SSDPs), where a convex subproblem/approximation of the original problem is solved at each iteration of the algorithm. The first issue with these schemes is that they often lack any good theoretical convergence guarantee (for instance, when the only guarantee is that the sequence of objective values is monotonously non-increasing). The second issue is that they may face significant practical difficulties, probably in terms of implementation, but more importantly in terms of numerical conditioning issues that may lead to a slow or no progress of the objective function.

Instead, what should have been considered are genuine non-convex optimization methods, i.e. methods that have a guarantee of convergence to stationary points and that are efficient in practice. There are two main classes of such methods: the gradient-based methods and the derivative-free methods. The current points are presented in a still summarized but more detailed fashion in \cite{SW12B}.

$\circ$ The trend of exploiting gradient-based methods to solve fundamental non-convex problems of optimization in systems and control is not new. In our opinion, the most significant such results have been developed in the 2000s in parallel by Michael Overton and Pierre Apkarian and their coworkers. Relying on Rockafellar \cite{RW98} and Clarke \cite{C83} non-smooth variational analysis, they developed non-smooth optimization methods that converge to locally optimal solutions for fixed-order (thus non-convex) control design problems. Two particularly visible contributions are the control design methods HIFOO \cite{BHLO06} (using the solver HANSO) and \texttt{hinfstruct} \cite{AN062}. The importance of these results and methods must be emphasized, and we strongly recommend the many related fundamental works of these authors.

$\circ$ Next to the gradient-based approach, there exist a second alternative that is apparently nearly unknown in systems and control. More precisely, what appears currently unacknowledged in this field is that there now exist strong convergence guarantees for a class of derivative-free optimization methods, thanks to which it is now justifiable to rely on such methods to solve non-convex optimization problems. For a one-paragraph summary of these guarantees, see \cite{SW12B}. For a fuller presentation, see the book \cite{CSV09}, the first ``\textit{comprehensive treatment}' of derivative-free methods. Note however that using these methods makes sense when the gradient information (of part of the objective function or the constraints) is not available (at all, or accurately, or at an acceptable computing cost). Also ideally, to keep short computing times with current desktop computers (e.g. $<$2-3 min), the number of variables should be reasonable (e.g. $<$25) and/or the objective function should be evaluated efficiently (e.g. $<$1s), or the method may be stopped before full convergence.

Several examples of the efficiency of derivative-free methods on control design problems are available in the (soon-to-be released, in September 2012) PhD thesis \cite{S12}. An appetizer is available in \cite{SW12B} (see hyperlinks and files therein), other different benchmarks (for which HIFOO and \texttt{hinfstruct} can be used) are given in \cite{S11} (see also the example in the hyperlinked presentation). Such examples will be drawn in the presentation, which will also be posted on arXiv.

To prompt a lively debate, we provide hereunder a 2 or 3 steps method to produce a virtually infinite amount of papers. This method relies on the assumption that obtaining a convex problem necessarily imply some measure of approximation of the original (engineering) problem (what is convex in this imperfect world?).
\begin{enumerate}
	\item Pick a paper where a problem is solved via LMI/convex optimization.
	\item (Optional) Rewrite the problem, so as to be less conservative with respect to the objective actually sought, and/or the constraints actually encountered.
	\item In the typical situation where the problem is non-convex, solve it via an adequate non-convex optimization method.
\end{enumerate}
The step 2. is optional in the sense that the problem in the step 1. may already be non-convex, but that a convex approach (such as an ILMI/SSDP) is proposed for it.

The convex approach really becomes unavoidable when the number of variables of the original problem is too great and/or when the actual problem is convex (which occurrences are often largely academic). We also remark that the solution returned by the convex optimization can be used as initial solution for the non-convex problem (as suggested e.g. in \cite[Subsec. 1.4.3]{BV04}). This may be helpful, or not. That is, randomly chosen (but feasible) initial solutions may be just as useful, in which case the convex optimization can be skipped entirely.

\end{document}